\documentclass[11pt]{amsart}
\usepackage{amsmath}
\usepackage{amsthm}
\usepackage{graphicx}

\newtheorem{theorem}{Theorem}
\newtheorem{proposition}[theorem]{Proposition}
\newtheorem{lemma}[theorem]{Lemma}
\newtheorem{corollary}[theorem]{Corollary}
\newtheorem{conjecture}[theorem]{Conjecture}
 
\theoremstyle{remark}

\newtheorem*{example}{Example}

\theoremstyle{definition}

\def\e{\varepsilon} 
\def\vol{\rm vol\mbox{ }} 
\def\qvol{{\rm vol}_q\mbox{ }} 
\def\1vol{{\rm vol}_{q=1}\mbox{ }} 
\def\0vol{{\rm vol}_{q=0}\mbox{ }} 

 \def\a{{\mathcal{A}}}
\def\cha{{\chi_\mathcal{A}}} 
 \def\ca{\widetilde{\chi}_\mathcal{A}}
\def\chb{{\chi_{\b}}} 
 
  \def\b{\mathcal{B}_n}
  \def\c{\mathbb{C}}
   \def\r{\mathbb{R}}
   \def\z{\mathbb{Z}}
\def\m{\mathcal{M}}
\def\base{\mathcal{B}}
\def\o{\mathcal{O}}

\def\ms{\mathcal{S}}

\def\i{\mathcal{I}}

\def\x{\mbox{\bf{x}}}
\def\s{\langle S \rangle}
\def\f{\varphi}

\begin{document}
\title[Branched polymers and hyperplane arrangements]{Branched polymers and hyperplane arrangements}
  
\author{Karola M\'esz\'aros  \hspace{.3in} Alexander Postnikov}
\address{
Department of Mathematics, Massachusetts Institute of Technology, Cambridge, MA 02139 
\newline
{\tt karola@math.mit.edu, apost@math.mit.edu}
}
\date{December 16, 2009}
\keywords{polymers, braid arrangement, hyperplane arrangement, characteristic polynomial, broken circuit, Orlik-Solomon algebra}

\begin{abstract}
 We generalize the construction of  connected  branched polymers and the notion of the volume of the space of connected branched polymers  studied by Brydges and Imbrie \cite{bi}, and Kenyon and Winkler \cite{kw} to any    hyperplane arrangement $\a$. The volume  
 of the resulting configuration space of connected branched  polymers  associated to the    hyperplane arrangement   $\a$    is expressed through the value of the characteristic polynomial of $\a$ at $0$.  
 We give a more general definition of the space of branched polymers, where we do not require connectivity, and introduce the notion of $q$-volume for it, which is expressed through the value of the characteristic polynomial of $\a$ at $-q$.   Finally, we relate the volume of the space of branched polymers to broken circuits and show that the cohomology ring of the space of branched polymers is isomorphic to   the Orlik-Solomon algebra.      \end{abstract}

\maketitle
  
  \section{Introduction}
  \label{sec:intro}
  
  Brydges and Imbrie \cite{bi}, and Kenyon and Winkler \cite{kw} study the space of  branched polymers of order $n$ in $\r^2$ and $\r^3$. In their definition, a  {\bf branched polymer} of order $n$ in $\r^D$ is a connected set of $n$ labeled unit spheres in $\r^D$ with nonoverlapping interiors. To each polymer $P$ corresponds a graph $G_P$ on the vertex set $[n]$ where $(i, j) \in E(G_P)$ if and only if spheres  $i$ and $j$ touch in the polymer. The space of  branched polymers of order $n$ in $\r^D$ can be parametrized by the $D$-dimensional angles at which the spheres are  attached to each other.  If $G_P$ has a cycle, this parametrization becomes ambiguous.  However, since polymers containing a cycle of touching spheres have probability zero, this ambiguity is not of concern.  Brydges and Imbrie \cite{bi}, and Kenyon and Winkler \cite{kw} compute the volume of the space of  branched polymers in $\r^2$ and $\r^3$, which are $(n-1)! (2 \pi)^{n-1}$ and $n^{n-1} (2 \pi)^{n-1}$, respectively.  They show that if we relax the requirement that the spheres in the polymer have radii $1$, the volume of the  space of   branched polymers in $\r^2$ remains the same, while the volume of the  space of   branched polymers in $\r^3$ changes. 
  
  Intrigued by the robustness of the space of branched polymers in $\r^2$ under the change of radii, we 
  generalize this notion differently from how it is done in  \cite{bi, kw}. Under our notion of polymers,  the volume of the  space of polymers is independent of the radii in all cases. 
  
  We associate the space of  branched polymers to any    hyperplane arrangement $\a$. The polymers corresponding to the braid arrangement $\b$  coincides with the definition of branched polymers in $\r^2$ given above.  We also broaden the notions of branched polymers, by  not requiring  polymers to be connected. We define the volume of the space of such polymers and show that it is invariant under the change of radii. In the case of the braid arrangement, the volume of the space of connected branched polymers associated to $\b$  is $(n-1)! (2 \pi)^{n-1}=(-2 \pi)^{r(\b)} {\chi}_{\b}(0)$, where $r(\b)$ is the rank of $\b$ and $\chb(t)$ is the characteristic polynomial of $\b$. In our generalized notion of volume, weighted with $q$, 
  the volume of the space of   branched polymers associated to $\b$ is
  $(-2 \pi)^{r(\b)} {\chi}_{\b}(-q)$.  The volume of the space of connected branched polymers associated to $\b$ is a specialization of this $q$-volume at $q=0$. 
 
  A theorem of the same flavor holds for any of our polymers.   
  \vspace{.05in}

\noindent {\bf Theorem.} {\it The $q$-volume of the space of branched polymers associated to a    central hyperplane arrangement $\a$ is 
   $(-2 \pi)^{r(\a)}{\chi}_{\a}(-q)$.   Furthermore, the volume of the space of connected branched polymers associated to $\a$ is a specialization of its $q$-volume at $q=0$. }
     \vspace{.05in}
     
   In the case that $\a$ is a graphical arrangement, we recover the $G$-polymers of \cite{kw}, and the Theorem  can be rephrased in terms of the chromatic polynomial of the graph $G$. We also relate the volume of  the space of branched polymers to broken circuits and the Orlik-Solomon algebra. Finally, we prove that the  cohomology ring of the space of branched polymers is isomorphic to  the Orlik-Solomon algebra, and  conjecture the same for for  the space of connected branched polymers.
  
  The outline of  the paper is as follows. In Section \ref{sec:braid} we explain how the notion of branched polymers   in $\r^2$ from \cite{bi, kw} translates to  polymers associated to the braid arrangement $\b$.     In Section \ref{sec:arr} we give a general definition of  connected branched polymers, as well as branched polymers (where we do not require connectivity) associated to any    arrangement $\a$. We define the notion of $q$-volume of the space of branched polymers and restate the Theorem  about the  value of the $q$-volume of the space of  branched polymers. In Section \ref{sec:proof} we prove the Theorem. In Section \ref{sec:gpoly} we recover the $G$-polymers of \cite{kw}   from the graphical arrangement $\a_G$. We relate the volumes of branched polymers to broken circuits in Section \ref{sec:nbc}. We conclude in Section \ref{sec:os} by proving that the cohomology ring of the space of  branched polymers is the Orlik-Solomon algebra, and we   conjecture that  the cohomology ring of  the space of connected branched polymers is the same.

  \section{Branched polymers and the braid arrangement}
\label{sec:braid}
 
 In this section we explain how to think of the space of  connected branched polymers and its volume defined in \cite{bi, kw} in terms of the braid arrangement $\b$. We  also give a more general definition  of the space of   branched polymers associated to $\b$ equipped with the notion of $q$-volume. For $q=0$ we recover the volume of the space of  connected  branched polymers in the sense of \cite{bi, kw}. \footnote{Note that our terminology differs slightly from that of \cite{bi, kw}. The notion of ``branched polymer" in \cite{bi, kw} corresponds to our ``connected branched polymers."}

Let $V:=\c^n / (1, \ldots, 1) \c \cong \c^{n-1}$. A point $\x \in V$ has coordinates $(x_1, \ldots, x_n)$ considered modulo simultaneous shifts of all $x_i$'s.

   The {\bf braid arrangement} $\b$  in $V$, $n >1$,   consists of the ${n \choose 2}$ hyperplanes   given by the equations $$h_{ij}(\mbox{\bf{x}})=x_i-x_j=0, \mbox{ for } 1\leq i<j\leq n, \mbox{ where } {\bf x}  \in V.$$ 
   
   Let  $$H_{ij}=\{\x \in V \mid h_{ij}(\x)=0\}.$$ 

Define   the {\bf space of branched polymers} 
 associated to the arrangement $\b$ and nonnegative scalars $R_{ij}$, $1\leq i<j\leq n$ to be 
   $$P_{\b}=\{ \x \in V  \mid  |h_{ij}(\x)|\geq R_{ij} \mbox{ for all } H_{ij} \in \b \}.$$  
 
 A connected branched polymer of size $n$ is a connected  collection of $n$ labeled disks in $\r^2$ with nonoverlapping interiors (considered up to translation).  Think of the collection of $n$ disks in $\r^2=\c$ as a point $\x \in V,$ where $x_k$ is the center of the $k^{th}$ disk. 
 Denote by $r_k$  the radius of the $k^{th}$ disk and let $R_{ij}=r_i+r_j$.  The condition that the disks do not overlap can be written as $|x_i-x_j|\geq R_{ij}$.  Disks $i$ and $j$ touch exactly if $|x_i-x_j|= R_{ij}$.  Thus the space  $P_{\b}$  consists of points corresponding to branched polymers, which are   not necessarily  connected.

\begin{figure}[htbp] 
\begin{center} 
\includegraphics[scale=.75]{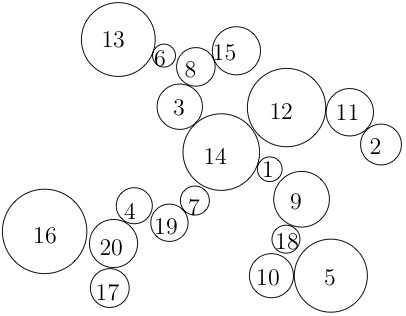} 
\caption{A connected branched polymer.} 
\label{fig0} 
\end{center} 
\end{figure}

 It is clear from the definition of $P_{\b}$ that it decomposes as the disjoint union 
 
 $$P_{\b}=\bigsqcup_{S \subset \b} P_{\b}^S,$$ where   
   $ P_{\b}^S$ consists of all points $\x \in P_{\b} $ such that  $|x_i-x_j|= R_{ij}$ exactly if  $H_{ij}  \in S$.  That is, 
  $$ P_{\b}^S=\{\x \in \c^n  \mid |h_{ij}(\x)|> R_{ij}\mbox{ for } H_{ij}\in \b \backslash  S,  \mbox{ }|h_{ij}(\x)|=R_{ij}\mbox{ for } H_{ij}  \in S  \}.$$

  We can consider  $\b$ as a   matroid with ground set the set of hyperplanes in $\b$ and independence given by linear independence of the normals to the hyperplanes. 
  
  Connected branched polymers, in the sense of \cite{bi, kw}, correspond to points in the strata $ P_{\b}^S$ of $P_{\b}$ such that $S$ contains a base of ${\b}$. Otherwise, the configuration of disks would not be connected. Thus,  the {\bf space of  connected branched polymers} corresponding to $\b$ is 
  
  $$CP_{\b}=\bigsqcup_{S \subset \b \mbox{ } : \mbox{ }  r(S)=r({\b})} P_{\b}^S,$$ where $r(S)$ denotes the rank of the arrangement $S$, which is the dimension of the space which the normals to the hyperplanes in $S$ span.

   We now define the notion of the volume of  the space $CP_{\b}$, which coincides with the definition given in \cite{bi, kw}.
   
   Let  $S=\{h_{i_1 j_1}, \ldots, h_{i_{n-1} j_{n-1}}\}$ be a base of $\b$.  Let us embed  $P_{\b}^S$ into $[0, 2\pi]^{n-1}$ by $$\x \mapsto \f(\x)= (\f_1, \ldots, \f_{n-1}),$$ where $h_{i_k j_k}(\x)=e^{i \f_{k}} R_{i_k j_k},  k \in [n-1] $. Define the volume of $P_{\b}^S$ as the volume of its image in $[0, 2\pi]^{n-1}$:

   $$\vol P_{\b}^S=\vol (\{ \f(\x) \mid \x \in  P_{\b}^S\}).$$  
    
     If $S$ is a dependent set in $\b$, then $P_{\b}^S$ has a lower dimension (${\rm dim} \mbox{ } P_{\b}^S<n-1$), so we let $\vol P_{\b}^S=0$.  
  
  Let
  $$\vol CP_{\b}=\sum_{S \mbox{ } : \mbox{ }  r(S)=r({\b})}\vol P_{\b}^S=\sum_{S \mbox{ } : \mbox{ }  S \in \base({\b})} \vol P_{\b}^S,$$ where $\base(\m)$ denotes the set of bases of a matroid $\m$.

   Recall that   the characteristic polynomial  of $\b$ is $\chb(t)=(t-1)\cdots (t-(n-1))$ \cite[p. 414]{s}.    
  By \cite[Theorem 2]{kw}, $\vol CP_{\b}=(n-1)! (2 \pi)^{n-1}$. Observe that it   equals    
  \begin{equation} \label{vol} \vol CP_{\b} =(-2 \pi)^{r(\b)}  \chb(0),\end{equation} where $r(\b)=n-1$ is the rank of $\b$.  Equation (\ref{vol}) is a special case of  Theorem \ref{01} which we state and prove in the following sections.

  \section{Polymers associated to a hyperplane    arrangement}
\label{sec:arr}
 
 In this section we associate branched polymers to any   central hyperplane  arrangement. We calculate the volume of the space of  connected branched polymers as well as the $q$-volume of the space of all branched polymers. 
  
 Let $\{h_1, \ldots, h_N\} \subset (\c^r)^*$, 
 and assume that $h_1, \ldots, h_N$ span $(\c^r)^*$. Let $\a=\{H_1, \ldots, H_N\} $ be the hyperplane arrangement where the hyperplanes are $$H_i=\{ \x \in \c^r \mid h_{i}(\x)=0\}.$$ Note that $r=r(\a)$ is the rank of the arrangement.
   
  Define the    {\bf space of branched polymers  } associated to the arrangement $\a$ and nonnegative scalars $R_{i}$, $1\leq i\leq N$ to be 
  $$P_{\a}=\{ \x \in \c^r  \mid |h_{i}(\x)|\geq R_{i} \mbox{ for all } H_i \in \a \}.$$    This space can be thought of as $\c^r$ with $N$ tubes removed. Namely, let $$T_i=\{ \x \in \c^r \mid |h_i(\x)|<R_i\}, \mbox{ } i \in [N],$$ be the $i^{th}$ tube, which is the set of points in $\c^r$ at distance less than $R_i/ ||h_i||$ from the hyperplane $H_i$. Clearly then $P_{\a}$ is the complement to all the tubes in $\c^r$.
The space $P_\a$ is related to the well-studied space $C_\a=\c^r \backslash \bigcup_{H \in \a} H,$  the complement of the arrangement $\a$ in $\c^r$.

 The space of branched polymers $P_{\a}$ decomposes as the disjoint union 
 
 \begin{equation} \label{bp} P_{\a}=\bigsqcup_{S \subset \a} P_{\a}^S,\end{equation} where 
  
  $$ P_{\a}^S=\{\x \in \c^r  \mid |h_{i}(\x)|> R_ {i}\mbox{ for } H_i  \in \a \backslash S, \mbox{ }  |h_{i}(\x)|=R_ {i} \mbox{ for } H_i \in S  \}.$$
  
  Consider  $\a$ as the matroid with ground set the  hyperplanes in $\a$ and independence linear independence of the normals to the hyperplanes. 
  
 Define the {\bf space of connected  branched polymers} associated to $\a$ as \begin{equation} \label{conn} CP_\a=\bigsqcup_{S \subset \a \mbox{ } : \mbox{ } r(S)=r(\a)} P_{\a}^S. \end{equation}  
 \medskip
 
  We  now   define the notions of volume  $\vol$and  $q$-volume $\qvol$of   $ P_{\a}^S$ for any set $S \subset \a$ with respect to $\a$. Let $\vol P_{\a}^S=0$ for any dependent set $S$, since ${\rm dim } \mbox{ } P_{\a}^S< n$ . 
  
   Let  $S=\{h_{i_1 }, \ldots, h_{i_{r} }\}$ be a base of $\a$.  Embed  $P_{\a}^S$ into $[0, 2\pi]^{r}$ by $$\x \mapsto \f(\x)=(\f_1, \ldots, \f_{r}),$$ where $h_{i_k }(\x)=e^{i \f_{k}} R_ {i_k },  k \in [r] $. Define the volume of $P_{\a}^S$ as the volume of its image in $[0, 2\pi]^{r}$:

   $$\vol P_{\a}^S=\vol (\{ \f(\x) \mid \x \in  P_{\a}^S\}).$$

  For any subset $S=\{H_{i_1}, \ldots, H_{i_l}\} $   of $\a$, let $$\langle S \rangle_{\a}:=\{ H_i \in \a \mid h_i \mbox{ is in the span } \langle h_{i_1}, \ldots, h_{i_l} \rangle \}.$$ We consider $\langle S \rangle_{\a}$
  as a hyperplane arrangement in $\c^r /  (H_{i_1} \cap \ldots \cap H_{i_l}) \cong \c^{r'}.$ Its rank is  $r(\langle S \rangle_{\a})=r'$.

Define

  $$\vol P_{\a}^S:=\vol P_{\langle S \rangle_{\a}}^S \cdot (2\pi)^{r(\a)-r(\langle S \rangle_{\a})},$$ and

  $$\qvol P_{\a}^S:=\vol P_{\a}^S \cdot q^{r(\a)-r(\langle S \rangle_{\a})}.$$

  Finally, define
    $$\vol CP_{\a}=\sum_{S \in \base(\a)}  \vol P_{\a}^S,  \mbox{  }\vol P_{\a}=\sum_{S \in \i(\a)} \vol P_{\a}^S, $$ and 
  $$\qvol P_{\a}=\sum_{S \in \i(\a)} \qvol P_{\a}^S, $$ where $\i(\m)$ denotes the independent sets of the matroid $\m$.
  \medskip

\begin{theorem} \label{01}
  \begin{equation} \label{beta1}
   \qvol P_\a=(-2 \pi)^{r(\a)} {\chi}_{\a}(-q),\end{equation}
  
 where $\cha(t)$ is the characteristic polynomial of the  central   hyperplane arrangement $\a$ in $\c^r$ with $\dim(\bigcap_{H \in \a} H)=0$. 
   
   {\it In particular,} 
   \begin{equation} \label{nobeta}
   \vol CP_\a=(-2 \pi)^{r(\a)} {\chi}_{\a}(0).\end{equation}
   \end{theorem}
   \medskip
   We  prove Theorem \ref{01} in the next section.
 
 \medskip
 
 We now point out an alternative way of  thinking about the $q$-volume. 
   Recall  that $\s_\a$ has the structure of a  matroid. Two hyperplanes  $H_k$ and $H_l$ are in the same   {\bf connected component} of $\s_\a$ if they are  in a common circuit of the matroid defined by $\s_\a$. Let the components of $\s_\a$ be $C_1, \ldots, C_k$. Then it is possible to write $S=\bigsqcup_{l=1}^k S_l$ such that $C_i=\langle S_i \rangle_\a$. We call $S_1, \ldots, S_k$ the {\bf components} of $S$.  
 
 The volume $\qvol P_{\langle S \rangle_{\a}}^S$, where $S_1, \ldots, S_k$ are the  components of $S$, can be written as 
 
 $$\qvol P_{\langle S \rangle_{\a}}^S=\qvol P_{\langle S_1 \rangle_\a}^{S_1} \cdots \qvol P_{\langle S_k \rangle_\a}^{S_k},$$ where $S_1, \ldots, S_k$ are the  components of $S$.

 \section{The proof of Theorem \ref{01}}
 \label{sec:proof}
 
 In this section we prove Theorem \ref{01}. We first state a few lemmas we use in the proof.

  \begin{lemma} \label{inv}  (Invariance Lemma) The volume  $  \sum_{S \mbox{ }: \mbox{ } |S|=s,  \mbox{ } S \in \i({\a})}  \vol P_{\a}^S$ is independent of the $R_ i$, $i \in [N]$, where $s \in \{0, 1, \ldots, r(\a)\}$.
 \end{lemma}
 
Lemma \ref{inv} is more general than the Invariance Lemma stated and proved in \cite{kw}, but the techniques used in \cite{kw} apply just as well in this case, yielding the proof of Lemma \ref{inv}.

 \begin{lemma} \label{a'} Let  $\a'=\a \backslash \{H_1\}$,   such that $r(\a)=r(\a')$. Then as  $R_ 1 \rightarrow 0$,  $$\vol P_{\a'}^{S}=\vol P_{\a}^{S},$$  for any  independent set $S\subset \a,$ with $H_1 \not \in S$. 
 
 \end{lemma} 
 
 \proof Note that $r(\s_\a)=r(\s_{\a'}).$ 
 Let the set of components of  $\s_{\a}$  be  $\{\langle S_1 \rangle_{\a}, \ldots, \langle S_k \rangle_{\a} \} $, and the set of components of $\s_{\a'}$  be  $\{ \langle S_{1i} \rangle_{\a'}\}_{i \in I_1} \cup \ldots \cup \{\langle S_{ki} \rangle_{\a'}\}_{i \in I_k}.$ Then $\bigcup_{i \in I_l} S_{li}=S_l$, $ l \in [k]$, since if two  elements $x$ and $y$ in the same component  of $\s_{\a'}$, then  there is a circuit in $\a$ containing $x, y$ and not containing $H_1$.  Also, if for two elements $x$ and $y$ in different components  of $\s_{\a'}$, there is a circuit in $\a$ containing $x$ and $y$ (then this circuit necessarily contains $H_1$) and if  $x'$ and $x$ in a circuit of $\a$ and $x$ and $y$ are in a circuit of $\a$, so are $x'$ and $y$. 
 
 We have   $$\vol P_{\a}^{S}=\prod_{j=1}^k  \vol P_{\langle S_j \rangle_\a}^{S_j} \cdot (2 \pi)^{r(\a)-r(\s_\a)}$$  and       $$\vol P_{\a'}^{S} =\prod_{j=1}^k \prod_{i \in I_j} \vol P_{\langle S_{ji} \rangle_{\a'}}^{S_{ji}}\cdot (2 \pi)^{r(\a')-r(\s_{\a'})}.$$

We claim that for any $j \in [k]$,   \begin{equation} \label{equal} \vol P_{\langle S_j \rangle_{\a}}^{S_j}=\prod_{i \in I_j} \vol P_{\langle S_{ji} \rangle_{\a'}}^{S_{ji}}.\end{equation} Since $r(\a)-r(\s_\a)=r(\a')-r(\s_{\a'})$, equation (\ref{equal})  suffices to prove  $\vol P_{\a}^{S}=\vol P_{\a'}^{S}$. 

Equation (\ref{equal})  follows, since, as  $R_ 1 \rightarrow 0$, the presence of the hyperplane $H_1$ in $\a$ imposes no constraints on the angles $\f_{j}$, for $j \in S_{ji}$. Furthermore, if $ H_l \in \langle S_{jl} \rangle_{\a'}$ and $ H_m \in \langle S_{jm} \rangle_{\a'}$, $l \neq m$,  then as $R_ 1 \rightarrow 0$, the span of the angles $\f_l$ and $\f_m$ are independent of each other as all circuits containing both of them also contain $H_1$.  \qed
 
  \begin{lemma} \label{a''} Let  $\a''=\a / \{H_1\}$ be the contraction of $\a$ with respect to the hyperplane $H_1$   such that $r(\a)=r(\a'')+1$. Then as  $R_ 1 \rightarrow 0$,  $$2 \pi \cdot \vol P_{\a''}^{S\backslash \{H_1\}}=\vol P_{\a}^{S},$$  for any independent set $S\subset \a$ with $H_1 \in S$.
 
 \end{lemma} 

 \proof Note that since $S$ is an independent set in $\a$, so is $S\backslash \{H_1\}$ in $\a''$. Also, 
 $r(\s_\a)=r(\langle S\backslash \{H_1\} \rangle_{\a''})+1$ and  $r(\a)-r(\s_\a)=r(\a'')-r(\langle S\backslash \{H_1\} \rangle_{\a''}).$ Thus, it suffices to prove that 
 $2 \pi \cdot \vol P_{\langle S\backslash \{H_1\} \rangle_{\a''}}^{S\backslash \{H_1\}}=\vol P_{\langle S \rangle_{\a}}^{S}.$  
 
 Let the set of components of  $\s_{\a}$  be  $\{\langle S_1 \rangle_{\a}, \ldots, \langle S_k \rangle_{\a} \} $ with $H_1 \in S_1$. Note that any circuit in $\s_\a$ not involving $H_1$ automatically carries over to a circuit in $\s_{\a''}$ (easy to see if we assume, without loss of generality, that $h_1(\x)=x_n$). 
 Also, any circuit in $\s_{\a''}$  automatically lifts  to a circuit in $\s_{\a}$. Thus, the set of components of  $\langle S \backslash \{H_1\}\rangle_{\a''}$  is  $\{\langle S_{1j} \rangle_{\a''}\}_{j \in I_1} \cup \{\langle S_2 \rangle_{\a''}, \ldots,  \langle S_k \rangle_{\a''} \} $, where $\bigcup_{j \in I_1} S_{1j}=S_1$. As $R_ 1 \rightarrow 0$, the volumes $\vol P_{\langle S_i \rangle_{\a''}}^{S_i}=\vol P_{\langle S_i \rangle_{\a}}^{S_i},$ $i \in \{2, \ldots, k\}$.  Finally,  
 if $ H_l \in \langle S_{1l} \rangle_{\a
 }$ and $ H_m \in \langle S_{1m} \rangle_{\a}$, $l \neq m$,  then as $R_ 1 \rightarrow 0$, the span of the angles $\f_l$ and $\f_m$ are independent of each other as all circuits containing both of them also contain $H_1$. Thus, $2\pi \cdot \prod_{j \in I_1} \vol P_{\langle S_{1j} \rangle_{\a''}}^{S_1j}=\vol P_{\langle S_1 \rangle_{\a}}^{S_1},$
  the $2 \pi$ accounting for the fact that $H_1 \in S_1$ and as $R_ 1 \rightarrow 0$, $\f_1$ ranges between $0$ and $2\pi$ freely. 
\qed

   \medskip
 
 \noindent {\it Proof of Theorem \ref{01}.}   Equation (\ref{nobeta}) is a consequence of (\ref{beta1}) by the definitions of the spaces $P_\a$ and $CP_\a$ and their notions of volume. To  prove equation (\ref{beta1}) we consider two cases depending on whether there is a hyperplane $H \in \a$ such that $r(\a)=r(\a')$, for  $\a'=\a \backslash \{H\}$.  
  
 Suppose $\a$ is such that if we delete any hyperplane $H$ from it, obtaining the hyperplane arrangement $\a'=\a \backslash \{H\}$, then $r(\a')=r(\a)-1$. This means that the hyperplanes of $\a$ constitute a base. Then,  
  
   \begin{align*}\qvol  P_{\a}&=\sum_{S \mbox{ }: \mbox{ } S \in \i({\a})} q^{r(\a)-|S|} \vol P_{\a}^S\\ &= \sum_{i=0}^{r(\a)} { r(\a) \choose i} q^{r(\a)-i} (2 \pi)^{r(\a)} \\ &=(1+q)^{r(\a)} (2 \pi)^{r(\a)}\\  &=(-2 \pi)^{r(\a)} (-q-1)^{r(\a)} \\  &=(-2 \pi)^{r(\a)}\cha(-q) . 
   \end{align*} 
 
 The second equality holds since the independent sets of size $i$ are exactly the  $i$-subset of the base $\a$, of which there are ${ r(\a) \choose i}$  and $ \vol P_{\a}^S=(2 \pi)^{r(\a)}$ for any set $S$, since $\a$ is a base itself. The last equality holds since if $\a$ itself is a base, then $\cha(t)=(t-1)^{r(\a)}$.

 Suppose now that there exists a hyperplane $H \in \a$ such that $r(\a)=r(\a')$ for $\a'=\a \backslash \{H\}$. Assume without loss of generality that $H=H_1$. 
If $\a'$ and $\a''$ are the deletion and restriction of $\a$ with respect to the hyperplane $h_1(\x)=0$, then  $r(\a)=r(\a')=r(\a'')+1$. Recall that \begin{equation} \label{rec1} \cha(t)=\chi_{\a'}(t)-\chi_{\a''}(t).\end{equation}  
 
 \medskip

 Let $$R(\a, q)= (-2 \pi)^{r(\a)} {\chi}_{\a}(-q)$$ denote the right hand side of equation (\ref{beta1}).  By (\ref{rec1}) 
 
 $$R(\a, q)=R(\a', q)+2 \pi R(\a'', q).$$ In particular, 
 
 $$[q^p]R(\a, q)=[q^p]R(\a', q)+2 \pi [q^p]R(\a'', q),$$ for any integer power $p$ of $q$.
 
 To prove the theorem it suffices to show that the same recurrence holds for \begin{equation} \label{rec2} L(\a, q)=\sum_{S \mbox{ }: \mbox{ } S \in \i({\a})} q^{r(\a)-|S|} \vol P_{\a}^S,\end{equation} since the case when we cannot use the recurrence is exactly when removing any hyperplane from $\a$ reduces the rank, and we proved the validity of Theorem \ref{01} for it above.

Recurrence (\ref{rec2}) holds, since     \begin{align*}  \sum_{S \mbox{ }: \mbox{ } |S|=s,  \mbox{ } S \in \i({\a})} & \vol P_{\a}^S=  \sum_{S' \mbox{ }: \mbox{ }  |S'|=s, \mbox{ }   {S'} \in \i({\a'})}  \vol P_{\a'}^{S'} + \\ &+2 \pi \sum_{S'' \mbox{ }: \mbox{ }  |S''|=s-1,  \mbox{ }  {S''} \in \i({\a''})}  \vol P_{\a''}^{S''},\end{align*}
 since Lemma \ref{a'} and \ref{a''} imply
   
 \begin{equation} \label{notin} \sum_{S' \mbox{ }: \mbox{ }  |S'|=s, \mbox{ }   {S'} \in \i({\a'})}  \vol P_{\a'}^{S'} =\sum_{S \mbox{ }: \mbox{ } H_1 \not \in S, \mbox{ }  |S|=s,  \mbox{ } S \in \i({\a})}  \vol P_{\a}^S\end{equation} and 
 
 \begin{equation} \label{in}  2 \pi \sum_{S'' \mbox{ }: \mbox{ }  |S''|=s-1,  \mbox{ }  {S''} \in \i({\a''})}  \vol P_{\a''}^{S''}=\sum_{S \mbox{ }: \mbox{ } H_1  \in S, \mbox{ }  |S|=s,  \mbox{ } S \in \i({\a})}  \vol P_{\a}^S. \end{equation}
 \qed

 \begin{example} By Theorem \ref{01}  the   $q$-volume of the space of branched polymers associated to the braid arrangement is   \begin{align} \label{uh} \qvol \mbox{ } P_{\b}&=(-2 \pi)^{r(\b)} {\chi}_{\b}(-q)\\  \nonumber&=(2 \pi)^{n-1} (q+1)\cdots (q+(n-1)).  \end{align}  A special case of (\ref{uh}) is equation (\ref{vol}). 
  \end{example}
 
  \section{G-polymers}
  \label{sec:gpoly}
  
  In this section we discuss $G$-polymers, as defined in \cite[Section 3]{kw}, and show how they correspond to graphical arrangements in our setup. We also rewrite the $q$-volume of the generalized $G$-polymers in the language of  its chromatic polynomial, which is a special case of its Tutte polynomial. 
  
  Given a graph $G=([n], E)$, let $R_ {ij}$ be positive scalars for each edge $(i, j)$ of $G$. A  {\bf $G$-polymer}, following Kenyon and Winkler, is a configuration of points $(x_1, \ldots, x_n) \in \c^n$  such that 
  
  \begin{itemize}

  \item $|x_i-x_j|\geq R_ {ij}$, for any edge $(i, j)$ of $G$
 
 \item $ |x_i-x_j|= R_ {ij}$ for all edges $(i, j)$ of some  connected subgraph of $G$.
 \end{itemize}
  
    If $G$ is not connected, there are no $G$-polymers.
 The volume of the space of $G$-polymers is defined by the angles made by the vectors from $x_i$ to $x_j$, where $(i, j) \in E$ is such that $|x_i -x_j|=R_ {ij}$. See \cite{kw} for further details.
 
 Recall that the {\bf graphical arrangement} $\a_G$ corresponding to the graph  $G=([n], E)$     consists of the $|E|$ hyperplanes $$H_{ij}=\{\x \in \c^n \mid x_i-x_j=0\} \mbox{ for } (i, j) \in E, i<j.$$ We consider $\a_G$
  as a hyperplane arrangement in $\c^n /  ( \bigcap_{(i, j) \in E} H_{ij}).$  
 The space connected branched polymers associated to $\a_G$ as defined in (\ref{conn}) coincide with the space of $G$-polymers. 
 
Let  $\chi_G(t)$ be the chromatic polynomial  of  the graph $G$,  $\widetilde{\chi}_{G}(t)=\chi_G(t)/t^{k(G)}$, where $k(G)$ is the number of components of $G$ and $r(G)= n-k(G)$.
 A special case of  Theorem \ref{01} is the following proposition.
 
 \begin{proposition} \label{chrom} {\rm (\cite[Theorem 4]{kw})}
  \begin{equation} \label{ag} \vol CP_{\a_G}=(-2 \pi)^{r(G)} \widetilde{\chi}_{G}(0).\end{equation}
 \end{proposition}

Proposition \ref{chrom} follows from Theorem \ref{01} since  $\widetilde{\chi}_G(t)=\chi_{\a_G}(t)$ (\cite[Theorem 2.7.]{s}) and  $r(G)=r(\a_G)$.

 Recall also that the chromatic polynomial  $\chi_G(t)$ of a graph $G$ is a special case of its Tutte polynomial $T_G(x, y)$: $\chi_G(t)=(-1)^{r(G)} t^{k(G)}T_G(1-t, 0)$ (\cite[Theorem 6, p. 346]{b}). In this light, Proposition \ref{chrom} coincides with   \cite[Theorem 4]{kw}, which is stated in terms of the Tutte polynomial of $G$.

 The {\bf space of generalized $G$-polymers} is $P_{\a_G}$ as defined in (\ref{bp}).
A special case of  Theorem \ref{01} is the following proposition. 

\begin{proposition} \label{G}
 \begin{align*}\qvol P_{\a_G}&=(-2 \pi)^{r(G)} \widetilde{\chi}_{G}(-q).  \end{align*} 
 \end{proposition}
 
 As Proposition \ref{G} illustrates, the notion of $q$-volume of the space of generalized $G$-polymers allows for  recovering the entire chromatic polynomial of~ $G$.

\section{Symmetric polymers}
\label{sec:sym}

In this section we specialize our results to the type $B_n$ Coxeter arrangement. This arrangement naturally gives rise to polymers symmetric with respect to the origin. 

The {\bf type $B_n$ Coxeter arrangement} $\a(B_n)$, $n >1$,  in $\c^n$,  consists of the $2{n \choose 2}+n$ hyperplanes  $$H_{ij}^-=\{ \x \in \c^n \mid x_i-x_j=0\}, H_{ij}^+=\{ \x \in \c^n \mid x_i+x_j=0\},$$ $$H_{k}=\{ \x \in \c^n \mid x_k=0\}, \mbox{ for } 1\leq i<j\leq n, k \in [n].$$
We consider $\a(B_n)$
  as a hyperplane arrangement in $V=\c^n /  ( \bigcap_{H \in \a(B_n)} H).$

\medskip
The space of branched polymers associated to $\a(B_n)$ and nonnegative scalars $R_ {ij}^+, R_ {ij}^-, R_ k$, for $1\leq i<j\leq n, k \in [n] $ is, by definition, 

\begin{align*}P_{\a(B_n)}=\{ \x \in V  \mid |x_i+x_j|\geq R_ {ij}^+, |x_i-x_j|\geq R_ {ij}^-, |x_k|\geq R_ {k}, \\\mbox{ for } 1\leq i<j\leq n, k \in [n] \}.\end{align*}

To obtain symmetric polymers, place a disk of radius $r_k$ around the points $x_k$ and $-x_k$ in the plane for all $k \in [n]$. Let $R_{ij}^+=R_{ij}^-=r_i+r_j$ and $R_k=r_k$, for $1\leq i<j\leq n, k \in [n]$. Then the condition  
$|x_i+x_j|\geq R_{ij}^+$ ensures that the disks around $x_i$ and $-x_j$, and  $-x_i$ and $x_j$  do not overlap, the condition  
$|x_i-x_j|\geq R_ {ij}^-$ ensures that the disks around $x_i$ and $x_j$, and $-x_i$ and $-x_j$ do not overlap, and the condition  
$|x_k|\geq R_ {k}$ ensures that the disks around $x_k$ and $-x_k$   do not overlap.

Since the characteristic polynomial of $\a(B_n)$ is given by $$\chi_{\a(B_n)} (t)=(t-1)\cdot (t-3) \cdots (t-(2n-1)), \mbox{ \cite[p. 451]{s}}$$  Theorem 1 specializes to the following corollary.

\begin{corollary} \label{cor:sym}
 \begin{equation} 
 \qvol P_{\a(B_n)}=(2 \pi)^n(q+1)\cdot (q+3) \cdots (q+(2n-1)).\end{equation}

    In particular, 
  \begin{equation} 
   \vol CP_{\a(B_n)}=(2 \pi)^n 1\cdot 3 \cdots (2n-1).\end{equation}
   \end{corollary}
\medskip

\begin{figure}[htbp] 
\begin{center} 
\includegraphics[scale=.75]{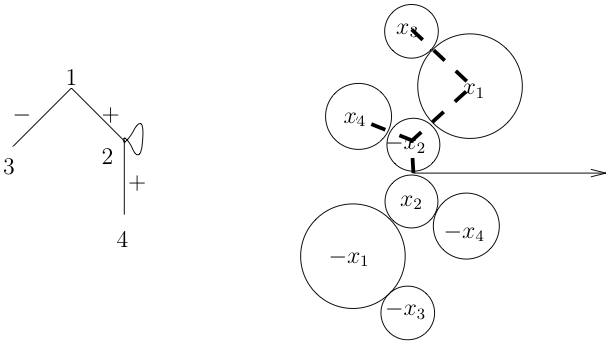} 
\caption{The graph representing a base and the corresponding polymer. The half-line drawn on the figure is the positive $x$-axis. The parametrization can be done by the angles which the dashed line segments make with the $x$-axis.} 
\label{fig1} 
\end{center} 
\end{figure}

\begin{figure}[htbp] 
\begin{center} 
\includegraphics[scale=.75]{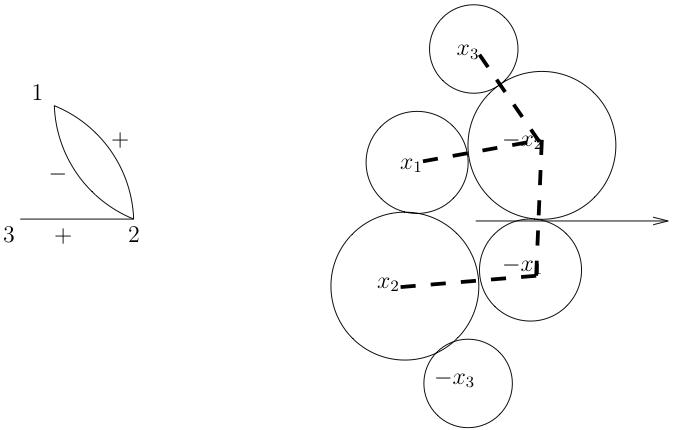} 
\caption{The graph representing a base and the corresponding polymer. The half-line drawn on the figure is the positive $x$-axis. The parametrization can be done by the angles which the dashed line segments make with the $x$-axis.} 
\label{fig2} 
\end{center} 
\end{figure}

We now outline the  geometric picture of connected branched polymers associated to $\a(B_n)$. 
Since $\vol CP_{\a(B_n)}=\sum_{S \in \base({\a(B_n)})} \vol P_{\a(B_n)}^S$ we restrict our attention to the polymers corresponding to the bases of $\a(B_n)$. Biject a subset $S$ of hyperplanes of $\a(B_n)$ with a graph $G_S$ on vertex set $[n]$  by adding an edge labeled by $+$ between vertices $i$ and $j$ if  $H_{ij}^+ \in S$,   by adding an edge labeled by $-$ between vertices $i$ and $j$ if  $H_{ij}^- \in S$, and adding a loop at vertex $k$ if $H_k \in S$. Then the bases of $\a(B_n)$ correspond to graphs $G$  on vertex set  $[n]$ with $n$ edges  such that each component of $G$ contains exactly one cycle or one loop, and if it contains a cycle then the cycle has an odd number of edges labeled by $+$. 

It is easy to see that if there is more than one loop in $G_S$, for some base $S$,  then the stratum $P_{\a(B_n)}^S$ is empty. This is true, since a loop based at vertex $k$  in $G_S$ corresponds to two disks of radius $R_ k=r_k$ around the points $x_k$ and $-x_k$, touching at the origin. If there was also a loop at vertex $l \neq k$, it would be impossible to also satisfy $|x_k+x_l|\geq r_l+r_k$ and $|x_k-x_l|\geq r_l+r_k$ if $k<l$ or $|x_l-x_k|\geq r_l+r_k$ if $k>l$. The disks corresponding to a component of $G_S$ with a loop at $k$ look like a tree symmetric across the origin and can be parametrized by the angles which the edges of  ``half" of the tree make with the $x$-axis and by the angle made by the segment $0-x_k$ and the $x$-axis; see Figure \ref{fig1}. The disks corresponding to a component of $G_S$ with a  cycle is symmetric with respect to the origin and contains a cycle of disks symmetrically around the origin. If the cycle is of length $2m$, then the polymer  is parametrized by $m+1$  angles which the edges of  the cycle make with the $x$-axis as well  the angles that  ``half" of the remaining edges  make with the $x$-axis; see Figure \ref{fig2}.

\section{Volumes of branched polymers and broken circuits}
\label{sec:nbc}

In this section we relate the volumes of branched polymers associated to $\a$ and broken circuits of $\a$. While the volume  $\vol \mbox{ }  CP_{\a}$ is invariant under changing the radii $R_i$ by the Invariance Lemma, the individual volumes of the strata $ P_{\a}^S$ change. 
By picking the radii $R_ i$ appropriately, it is possible to construct a stratification, where only strata corresponding to certain special bases $S$ appear in the stratification, and they all have the same volume. The mentioned special bases are the no broken circuit bases as we show in Theorem \ref{bc}.  

A {\bf broken circuit} of $\a$ with respect to a linear ordering $\o$ of the hyperplanes of $\a$ is a set $C-\{u\},$ where $C$ is a circuit and $u$ is the largest element of $C$ in the linear ordering $\o$. While the exact set of broken circuits depends on the linear order we choose, remarkably, the number of sets $S$ of hyperplanes  of size $k$, for an arbitrary fixed $k$, such that $S$ contains no broken circuit is the same; see \cite[Theorem 4.12]{s}.  

\begin{theorem} \label{bc}

Let  $\a=\{H_1, \ldots, H_N\}$ be a hyperplane arrangement with $R_ 1 \ll \cdots \ll R_N$. Then the space $CP_\a$ is a disjoint union of $(-1)^{r(\a)}\cha(0)$ tori $(S^1)^{r(\a)}$.

More precisely, if $\o$ is the  order $H_1 <\cdots<H_N$  on the hyperplanes of $\a$ and $BC_\o(\a)$ are all subsets of $\a$ containing no broken circuits with respect to $\o$, then  for $R_ 1 \ll \cdots \ll R_N$

\begin{equation} \label{eq:tori} CP_{\a}=\bigsqcup_{S \in \base({\a}) \cap BC_\o(\a)} P_{\a}^S,\end{equation} where    $\vol \mbox{ }  P_{\a}^S=(2\pi)^{r(\a)}$ for $S \in \base({\a}) \cap BC_\o(\a)$.



 \end{theorem}

\proof  By  \cite[Theorem 4.12]{s}  

\begin{equation} \label{412} \# \{S   \mid |S|=i, S \in BC_\o(\a) \} =(-1)^{r(\a)}[q^{r(\a)-i)}]\cha(-q).\end{equation} Thus, $|\base({\a}) \cap BC_\o(\a)|=(-1)^{r(\a)}\cha(0)$, and it suffices to prove the second
 part of the theorem's statement.

Recall that by definition

$$ CP_\a=\bigsqcup_{S \subset \a \mbox{ } : \mbox{ } r(S)=r(\a)} P_{\a}^S. $$

We first show that if $S$ contains a broken circuit, then $P_{\a}^S$ is empty. Let $\{H_{i_1}, \ldots, H_{i_k}\}  \subset S$, $i_1<\cdots <i_k$,  be a broken circuit. Then there exists   $H_i \in \a$  such that $\{H_{i_1}, \ldots, H_{i_k}, H_i\}$,    is a circuit and $i_k<i$. Let  $h_i =\sum_{j=1}^k c_j h_{i_j}$,   $c_j\neq 0$.  Then for any $\x \in P_{\a}^{S_j}$, $$\sum_{j=1}^k |c_j| R_{i_k} \geq \sum_{j=1}^k |c_j| R_{i_j}=\sum_{j=1}^k |c_j h_{i_j}(\x)|\geq |\sum_{j=1}^k c_j h_{i_j}(\x)|=|h_i(\x)|\geq R_{i}.$$ However, $$\sum_{j=1}^k |c_j| R_{i_k} < R_{i}, \mbox{  since } R_{i_k} \ll R_i.$$ Thus, $P_{\a}^{S}$ is empty.

Next we prove that  $\vol P_{\a}^S=(2\pi)^{r(\a)}$ for $S \in \base({\a}) \cap BC_\o(\a)$. Let $S=\{H_{j_1}, \ldots, H_{j_r}\}$ with ${j_1}< \cdots<{j_r}$. Since $S$ contains no broken circuit, it follows that if $ i \in [N]\backslash\{{j_1}, \ldots, {j_r}\}$, then $i<j_r$. To prove that $\vol P_{\a}^S=(2\pi)^{r(\a)}$, it suffices to show that $|h_{j_l}(\x)|=R_{j_l}$ for $l \in [r]$ imply $|h_{i}(\x)|>R_{i}$ for $ i \in [N]\backslash\{{j_1}, \ldots, {j_r}\}$. Since $S$ is a base, $S\cup \{H_i\}$ contains a circuit $C$ whose maximal element is $H_{j_k}\neq H_i$,  $k \in [r]$, since $S \in   BC_\o(\a)$.
Thus, $$h_{j_k}(\x)=c_i h_i(\x)+\sum_{l \in [N]\backslash \{k\}}c_{j_l}h_{j_l}(\x),$$ where $c_i, c_{j_1}, \ldots, c_{j_r}$  are scalars  with $c_i \neq 0$.  In particular, $$|h_{j_k}(\x)| \leq |c_i|| h_i(\x)|+\sum_{l \in [N]\backslash \{k\}}|c_{j_l}||h_{j_l}(\x)|.$$ Since $|h_{j_l}(\x)|=R_{j_l}$ for $l \in [r]$ and $R_{j_k} \gg R_i$
it follows that $|h_{i}(\x)|> R_{i}.$
  \qed 
 
  \vspace{.05in}

\noindent {\it An alternative proof of Theorem \ref{01}.} By definition, 
\begin{equation} \label{um2} \qvol P_\a=\sum_{i=0}^{r(\a)} \sum_{S \in \ms_i} (q \cdot 2\pi)^{r(\a)-i} \vol CP_{\s_\a},\end{equation} where $\ms_i$ is a collection of independent sets $S$ of cardinality $i$ such that the arrangements $\s_\a$ for $S \in \ms_i$ run over all rank $i$ subarrangements of $\a$. By Theorem \ref{bc} the right hand side of equation (\ref{um2}) is 

\begin{align} \label{um3}  \sum_{i=0}^{r(\a)} \sum_{S \in \ms_i} (q \cdot 2\pi)^{r(\a)-i} \# \{B   \mid B \in \base({\s_\a}) \cap BC_\o(\s_\a) \} \cdot (2\pi)^i= \\  \label{um4} \sum_{i=0}^{r(\a)}  q^{r(\a)-i}  \cdot (2\pi)^{r(\a)}\# \{S  \mid |S|=i, S \in BC_\o(\a) \}.   \end{align}

Since by  \cite[Theorem 4.12]{s}  

\begin{equation} \label{4.12} \# \{S   \mid |S|=i, S \in BC_\o(\a) \} =(-1)^{r(\a)}[q^{r(\a)-i)}]\cha(-q),\end{equation} 
it follows from equations (\ref{um2}), (\ref{um3}), (\ref{um4}) and (\ref{4.12}) that $$   \qvol P_\a=(-2 \pi)^{r(\a)} {\chi}_{\a}(-q).$$ In particular,
  $$\vol CP_\a=(-2 \pi)^{r(\a)} {\chi}_{\a}(0).$$

\qed

Note the striking similarity between the formula for the $1/t$-volume of branched polymers $P_{\a}$,  \begin{equation} \label{eq:vol}  {\rm vol}_{1/t} \mbox{ } P_{\a}=(-2 \pi)^{r(\a)}  {\chi}_{\a}(-1/t)\end{equation} 

and the generating function 

\begin{equation} \label{eq:os} \sum_{i\geq 0} {\rm rank } \mbox{ } H^i(C_\a, \z)t^i=(-t)^{r(\a)} \cha(-1/t),\end{equation}  where $C_\a=\c^r \backslash \bigcup_{H \in \a} H,$ is the complement of the arrangement $\a$ in $\c^r$ and $H^k(C_\a, \z)$ is the $k^{th}$ graded piece of the cohomology ring $H^*(C_\a, \z)$; see \cite{os}. 
The right hand sides of (\ref{eq:vol}) and (\ref{eq:os})  are in fact identical, if we normalize appropriately. 
  This is no surprise, since we related  ${\rm vol}_{1/t} \mbox{ }  P_{\a}$ to no broken circuits in Theorem \ref{bc}, and since the cohomology ring $H^*(C_\a, \z)$ of $C_\a$ is isomorphic to the Orlik-Solomon algebra $A(\a)$ associated to the hyperplane arrangement $\a$, and the $i^{th}$ graded piece of $A(\a)$ has a basis in terms of no broken circuits.

    \section{Cohomology rings of the spaces $P_\a$ and $CP_\a$ and the Orlik-Solomon algebra}
    \label{sec:os}
 
 In this section we prove and conjecture about the cohomology rings of the spaces $P_\a$ and $CP_\a$.  The well-known Orlik-Solomon algebra is isomorphic to the cohomology ring  $H^*(C_\a, \z)$ of $C_\a$, where  $C_\a=\c^r \backslash \bigcup_{H \in \a} H,$ is the complement of the arrangement $\a$ in $\c^r$. We prove that the spaces $P_\a$ and $C_\a$ are homotopy equivalent, and we conjecture that the same is true under certain circumstances for the spaces $CP_\a$ and $C_\a$. 
   
 \begin{proposition} \label{f}
 The map $f: C_\a \rightarrow P_\a$ defined by
 
 $$f : (x_1, \ldots, x_r) \mapsto {\rm max}_{i \in [r]}(1, \frac{R_ i}{|h_i(\x)|}) \cdot (x_1, \ldots, x_r),$$ where $\a$ is a    hyperplane arrangements in $\c^r$,  $$C_\a=\c^r \backslash \bigcup_{H \in \a} H,$$ $$P_{\a}=\{ \x \in \c^r  \mid |h_{i}(\x)|\geq R_ {i} \}$$ is a deformation retraction. 

 \end{proposition}
 
 \proof It is straightforward to check that $f$  satisfies the following three conditions:  
 
 (i) continuous,
 
 (ii) $f |_{P_\a}=id$,
 
 (iii) $f$ is homotopic to the identity.   
 \qed

  \begin{corollary} \label{corf}
The spaces  $P_\a$ and $C_\a$ are homotopy equivalent. In particular, the cohomology rings of  the spaces $C_\a$ and $P_\a$, where $\a$ is a    hyperplane arrangements in $\c^r$,    satisfy $$H^*(P_\a, \z)\cong H^*(C_\a, \z).$$ 
 \end{corollary}

 In Corollary \ref{bc1} we proved that for $\a=\{H_1, \ldots, H_N\}$ and  $R_1\ll \cdots \ll R_N$ the space $CP_\a$ is a disjoint union of $|\cha(0)|$ tori $(S^1)^{r(\a)}$. In particular, $CP_\a$ is disconnected and cannot be homotopy equivalent to $C_\a$.  However, if the scalars $R_i$, $i \in [N]$, satisfy certain generalized triangle inequalities,   we conjecture that the spaces $CP_\a$ and $C_\a$ are homotopy equivalent, so in particular, the  space  $CP_\a$ is connected.
  
   Let $\a=\{H_1, \ldots, H_N\} $ be the hyperplane arrangement where the hyperplanes are $H_i=\{ \x \in \c^r \mid h_{i}(\x)=0\}.$ Define the subset $T_\a$ of $ \r_{>0}^N$ to consist of all points  ${\bf R}=(R_1, \ldots, R_N) \in \r_{>0}^N$ satisfying the inequalities $$ R_{i_0} \leq |c_1| R_{i_1}+\cdots + |c_k|R_{i_k},$$ whenever $$ h_{i_0}=c_1 h_{i_1}+\cdots + c_k h_{i_k}.$$ 
   
   \begin{example} The set $T_{\b}$ consists of points ${\bf R}=(R_{ij})_{1\leq i< j\leq n}$, where $R_{ij}$ denotes the coordinate corresponding to the hyperplane $H_{ij}=\{ \x \in V \mid x_{i}-x_j=0\},$ for which the $R_{ij}$'s satisfy the triangle inequalities  $R_{ij} \leq R_{ik}+R_{kj}$, for all $1\leq i<  k<j \leq n$ and  $R_{ij} \leq R_{ik}+R_{jk}$, for all $1\leq i< j<k \leq n$. In particular,  $(R_{ij}=r_i+r_j)_{1\leq i< j\leq n} \in T_{\b}$, where $r_i$, $i \in [n]$, are nonnegative scalars.
   
   \end{example}

  \begin{conjecture} \label{?} For $(R_1, \ldots, R_N) \in T_\a$ the spaces $CP_\a$ and $C_\a$ are homotopy equivalent. In particular,  $H^*(CP_\a, \z)\cong H^*(C_\a, \z).$ 
  \end{conjecture}
     
 A proof of  Conjecture \ref{?} would yield another explanation, why, while $C_\a$ has $2n$ real dimensions, its top cohomology is in dimension $n$, since $CP_\a$ has $n$ real dimensions.

\end{document}